\documentclass[10pt]{amsart}
\usepackage{amsmath,amsthm,amscd,amsfonts,amssymb,graphicx,paralist,color}
\usepackage{enumitem}
\usepackage[initials]{amsrefs}

\newcommand{\qp}{\mathbb{Q}_p}

\newcommand{\zp}{\mathbb{Z}_p}

\newtheorem{theorem}{Theorem}[section]
\newtheorem{definition}[theorem]{Definition}
\newtheorem{proposition}[theorem]{Proposition}
\newtheorem{corollary}[theorem]{Corollary}
\newtheorem{lemma}[theorem]{Lemma}

\theoremstyle{definition}

\numberwithin{equation}{section}
\begin{document}
\title{Some properties of differentiable $p$-adic functions}
	\author[Fern\'andez]{J. Fern\'andez-S\'anchez}
	\address[J. Fern\'andez-S\'anchez]{\mbox{}\newline \indent 
		Grupo de investigaci\'on ``Teor\'ia de c\'opulas y aplicaciones'', \newline \indent 
		Universidad de Almer\'ia, \newline \indent 
		Carretera de Sacramento s/n, \newline \indent
		04120 Almer\'ia (Spain).}
	\email{juanfernandez@ual.es}
	
\author[Maghsoudi]{S. Maghsoudi}
\address[S. Maghsoudi]{\mbox{}\newline \indent Department of Mathematics,\newline \indent  University of Zanjan,\newline \indent  Zanjan 45371-38791, Iran}
\email{s\_maghsodi@znu.ac.ir}

	\author[Rodr\'{i}guez]{D.L.~Rodr\'{i}guez-Vidanes}
	\address[D.L.~Rodr\'{i}guez-Vidanes]{\mbox{}\newline\indent Instituto de Matem\'atica Interdisciplinar (IMI), \newline \indent Departamento de An\'{a}lisis y Matem\'{a}tica Aplicada \newline \indent Facultad de Ciencias Matem\'{a}ticas \newline \indent Plaza de Ciencias 3 \newline \indent Universidad Complutense de Madrid \newline \indent	Madrid, 28040 (Spain).}
	\email{dl.rodriguez.vidanes@ucm.es}
	
	\author[Seoane]{J.B. Seoane--Sep\'ulveda}
	\address[J.B. Seoane--Sep\'ulveda]{\mbox{}\newline\indent Instituto de Matem\'atica Interdisciplinar (IMI), \newline\indent Departamento de An\'{a}lisis Matem\'{a}tico y Matem\'atica Aplicada,\newline\indent Facultad de Ciencias Matem\'aticas, \newline\indent Plaza de Ciencias 3, \newline\indent Universidad Complutense de Madrid,\newline\indent 28040 Madrid, Spain.}
	\email{jseoane@ucm.es}

\thanks{The third and fourth authors were supported by Grant PGC2018-097286-B-I00. The third author was also supported by the Spanish Ministry of Science, Innovation and Universities and the European Social Fund through a “Contrato Predoctoral para la Formación de Doctores, 2019” (PRE2019-089135). The second author was supported by the Iran National Science	Foundation (INSF) grant no. 99019850.}

\begin{abstract}
	In this paper, using the tools from the lineability theory, we  distinguish certain  subsets of $p$-adic differentiable functions. Specifically,  we show that the following sets of functions are large enough to contain an infinite dimensional algebraic structure: (i) continuously differentiable but not  strictly differentiable functions, (ii) strictly differentiable functions of order $r$ but not strictly differentiable of order $r+1$, (iii) strictly differentiable functions with zero derivative  that are not Lipschitzian of any order $\alpha >1$, (iv) differentiable functions with unbounded derivative, and (v) continuous functions that are differentiable on a full set with respect to the Haar measure but not differentiable on its complement having cardinality the continuum.
\end{abstract}
\maketitle
\section{Introduction and terminology}

In the non-Archimedean setting, at least two notions of differentiability have been defined: classical and strict derivative. Classical derivative has some unpleasant and strange behaviors, but it has been long known that the strict differentiability in the sense of Bourbaki is the hypothesis most useful to geometric applications, such as inverse theorem. Let us recall  definitions.

Let $\mathbb K$ be a valued field containing $\qp$ such that $\mathbb K$ is complete (as a metric space), and  $X$ be a nonempty subset of $\mathbb K$ without isolated points.
Let $f\colon X \to \mathbb K$ and $r$ be a natural number. Set
	$$
	\nabla^{r} X:=\begin{cases}
	X & \text{if } r=1,\\
	\left\{\left(x_1, x_2,\ldots, x_r\right)\in X^r \colon x_i\neq x_j \text{ if } i\neq j \right\} & \text{if } r\geq 2.
	\end{cases}
	$$
The $r$-th difference quotient $\Phi_rf: \nabla^{r+1} X \to \mathbb K$ of $f$, with $r\geq 0$, is recursively given by $\Phi_0f=f$ and, for $r\geq 1$, $(x_1, x_2,\ldots, x_{r+1})\in \nabla^{r+1} X$ by 
	$$
	\Phi_r f(x_1, x_2,\ldots, x_{r+1})=\frac{\Phi_{r-1}f(x_1, x_3,\ldots, x_{r+1})-\Phi_{r-1}f(x_2, x_3,\ldots, x_{r+1})}{x_1-x_2}.
	$$
Then a function  $f\colon X \to \mathbb K$ at a point $a\in X$  is said to be:
\begin{itemize}
	\item \textbf{differentiable}  if     $\lim_{x\to a}f(x)-f(a)/(x-a)$ exists;
	
	\item \textbf{strictly differentiable of order $r$} if  $\Phi_r f$ can be extended to a continuous function $\overline{\Phi}_r  f: X^{r+1} \to \mathbb K$. We then set $D_rf(a)=\overline{\Phi}_r  f(a, a,\ldots, a)$.
\end{itemize}

Our aim in this work is to study these notions through a recently coined approach--the lineability theory. Searching  for large algebraic structures in the sets of objects with a special  property,  we, in this approach, can get  deeper understanding of the behavior of  the objects under discussion.   In \cites{jms1,jms2}  authors studied lineability notions in the $p$-adic analysis; see
also \cites{preprint-June-2019,jms3, fmrs}. The study of lineability  can be traced back to   Levine
and Milman \cites{levinemilman1940} in 1940, and  V. I. Gurariy \cite{gur1} in 1966.  These works, among
others, motivated the introduction of the notion of   lineability  in 2005 \cite%
{AGS2005} (notion coined by V. I. Gurariy). Since then it has been a rapidly developing trend in mathematical analysis. There are extensive works on the classical lineability theory, see e.g.  \cites{ar,AGS2005,ar2,survey,bbf,bay}, whereas some recent topics can be found in \cites{TAMS2020,abms,bns,cgp,nak,CS2019,bo,egs,bbls,ms}. It is interesting to note that Mahler in \cite{Ma2} stated that:
\begin{quote}
``On the other hand, the behavior of continuous functions of a $p$-adic variable is quite distinct from that of real continuous functions, and many basic theorems of real analysis have no $p$-adic analogues. 

$\ldots$ there exist infinitely many linearly independent non-constant functions the derivative of which vanishes identically $\ldots$''. 
\end{quote}  

Before further going, let us recall three essential notions. 
Let $\kappa$ be  a cardinal number. We say that a subset $A$ of a vector space $V$ over a field $\mathbb{K}$ is 
\begin{itemize}

	\item \textbf{$\kappa$-lineable} in $V$   if there exists a vector space $M$ of dimension $\kappa$ and $%
	M\setminus \{0\} \subseteq A$;
\end{itemize}	 
	and following \cites{ar2, ba3}, if $V$ is
	contained in a (not necessarily unital) linear algebra,  then $A$ is called
\begin{itemize}[resume*]	
	\item \textbf{$\kappa$-algebrable} in $V$ if there exists an algebra $M$ such that $%
	M\setminus \{0\} \subseteq A$ and $M$ is a $\kappa$-dimensional vector
	space;
	
	\item \textbf{strongly $\kappa$-algebrable} in $V$ if there exists an $\alpha$-generated free algebra $M$ such that $M\setminus \{0\}\subseteq A$.
\end{itemize}
Note that if $V$ is also contained in a commutative algebra, then a set $B\subset V$ is a generating set of a free algebra contained in $A$ if and only if for any $n\in \mathbb N$ with $n\leq \text{card}(B)$ (where $\text{card}(B)$ denotes the cardinality of $B$), any nonzero polynomial $P$ in $n$ variables with coefficients in $\mathbb K$ and without free term, and any distinct $b_1,\ldots,b_n\in B$, we have $P(b_1,\ldots,b_n)\neq 0$.

Now we can give an outline of our work. In Section \ref{bac} we  recall some standard concepts and notations concerning
non-Archimedean   analysis. Then, in the section of Main results, we 
show, among other things, that: 
(i) the set of  functions $\qp\rightarrow \qp$ with continuous derivative that are not strictly differentiable   is $\mathfrak{c}$-lineable ($\mathfrak{c}$ denotes the cardinality of the continuum), 
(ii) the set of strictly differentiable functions $\qp\rightarrow \qp$ of order $r$ but not strictly differentiable of order $r+1$ is $\mathfrak c$-lineable,
(iii) the set of all strictly differentiable functions $\zp\rightarrow \mathbb K$ with zero derivative  that are not Lipschitzian of any order $\alpha>1$ is $\mathfrak{c}$-lineable and $1$-algebrable, 
(iv) the set of differentiable functions $\qp\rightarrow \qp$ which derivative is unbounded is strongly $\mathfrak{c}$-algebrable, 
(v) the set of continuous functions $\zp\rightarrow \qp$ that are differentiable with bounded derivative on a full set for any positive real-valued Haar measure on $\zp$ but not differentiable on its complement having cardinality $\mathfrak{c}$ is $\mathfrak c$-lineable.

\section{Preliminaries for $p$-adic analysis}
\label{bac}

We   summarize  some basic definitions of $p$-adic analysis (for a
more profound treatment of this topic we refer the interested reader to \cites{go,kato,va,sc}).

We shall use standard set-theoretical notation. As usual, $\mathbb{N}, 
\mathbb{N}_0, \mathbb{Z}, \mathbb{Q}$, $\mathbb{R}$ and $\mathbb P$ denote
the sets of all natural, natural numbers including zero, integer, rational, real, and prime numbers, respectively.
The restriction of a function $f$ to a set $A$ and the characteristic function of a set $A$ will be denoted
by $f\restriction A$ and $1_A$, respectively.

Frequently,  we use a theorem due to Fichtenholz-Kantorovich-Hausdorff \cites{fich,h}: For any set $X$ of infinite cardinality there exists a family $\mathcal{B}\subseteq 
\mathcal{P}(X)$ of cardinality $2^{\text{card}(X)}$ such that for any finite
sequences $B_1,\ldots , B_n\in \mathcal{B}$ and $\varepsilon_1, \ldots,
\varepsilon_n \in \{0,1\}$ we have $B_1^{\varepsilon_1}\cap \ldots \cap
B_n^{\varepsilon_n}\neq \emptyset$, where $B^1=B$ and $B^0=X\setminus B$.  
A family of subsets of $X$ that satisfy the latter condition is called a family of \textit{independent subsets} of $X$.
Here $\mathcal{P}(X)$ denotes the power set of $X$. In what follows we fix  $\mathcal{N}$, $\mathcal N_0$ and $\mathcal P$ for families of independent subsets of $\mathbb{N}$, $\mathbb N_0$ and $\mathbb P$, respectively,
such that $\text{card}(\mathcal N)=\text{card}(\mathcal N_0)=\text{card}(\mathcal P)=\mathfrak c$.

Now let us recall that given a field $\mathbb K$, an absolute value on $\mathbb K$ is a function 
	$$
	|\cdot|\colon \mathbb K\to [0,+\infty)
	$$ 
such that:
	\begin{itemize}
		\item $|x|=0$ if and only if $x=0$,
		
		\item $|x y|=|x||y|$, and
		
		\item $|x+y|\leq |x|+|y|$,
	\end{itemize}
for all $x,y\in \mathbb{K}$. 
The last condition is the so-called \textit{triangle inequality}.
Furthermore, if $(\mathbb K,|\cdot|)$ satisfies the condition $|x+y|\leq \max\{|x|, |y|\}$ (the \textit{strong triangle inequality}), then $(\mathbb K,|\cdot|)$ is called a non-Archimedean field.
Note that $(\mathbb K,|\cdot|)$ is a normed space since $\mathbb K$ is a vector space in itself.
For simplicity, we will denote for the rest of the paper $(\mathbb K,|\cdot|)$ by $\mathbb K$.
Clearly, $|1|=|-1|=1$ and, if $\mathbb K$ is a non-Archimedean field, then $\underbrace{|1+\cdots+1|}_{n \text{ times}}\leq 1$ 
for all $n\in \mathbb{N}$. An immediate consequence of the strong triangle
inequality is that $|x|\neq |y|$ implies $|x+y|=\max\{|x|, |y|\}$.
Notice that if $\mathbb K$ is a finite field, then the only possible absolute value that can be defined on $\mathbb K$ is the trivial absolute value, that is, $|x|=0$ if $x=0$, and $|x|=1$ otherwise.
Furthemore, given any field $\mathbb K$, the topology endowed by the trivial absolute value on $\mathbb K$ is the discrete topology.

Let us fix a prime number $p$ throughout this work. For any non-zero integer 
$n\neq 0$, let $\text{ord}_p (n)$ be the highest power of $p$ which divides $%
n$. Then we define $|n|_p=p^{-\text{ord}_p (n)}$, $|0|_p=0$ and $|\frac{n}{m}%
|_p=p^{-\text{ord}_p (n)+\text{ord}_p (m)}$, the $p$-adic absolute value. The completion of the field of
rationals, $\mathbb{Q}$, with respect to the $p$-adic absolute value
is called the field of $p$-adic numbers $\mathbb{Q}_p$. 
An important property of $p$-adic numbers is that each nonzero $x\in \qp$ can be represented as
	$$
	x=\sum_{n=m}^\infty a_n p^n,
	$$
where $m\in \mathbb Z$, $a_n\in \mathbb F_p$ (the finite field of $p$ elements) and $a_m\neq 0$.
If $x=0$, then $a_n=0$ for all $n\in \mathbb Z$.
The $p$-adic absolute value
satisfies the strong triangle inequality. Ostrowski's Theorem states that every nontrivial
absolute value on $\mathbb{Q}$ is equivalent (i.e., defines the same
topology) to an absolute value $|\cdot|_p$, where $p$ is a prime number, or
the usual absolute value (see \cite{go}). 

Let $a\in \mathbb{Q}_{p}$ and $r$ be a positive number. The set $%
B_r(a)=\{x\in \mathbb{Q}_{p}\colon |x-a|_{p}<r\}$ is called the open ball of
radius $r$ with center $a$, $\overline{B}_{r}(a)=\{x\in \mathbb{Q}%
_{p}\colon |x-a|_{p}\leq r\}$ the closed ball of radius $r$ with center $a$, and $%
S_r(a)=\{x\in \mathbb{Q}_{p}\colon |x-a|_{p}=r\}$ the sphere of radius $r$ and
center $a$. 
It is important to mention that $B_r(a)$, $\overline{B}_{r}(a)$ and $S_r(a)$ are clopen sets in $\qp$.
The closed unit ball
 $$
 \mathbb{Z}_{p}=\{x\in \mathbb{Q}_{p}\colon |x|_{p}\leq 1\}=\left\{x\in \mathbb{Q}_{p}\colon x=\sum_{i=k}^\infty a_ip^i,\ a_i\in \{0,1,\ldots ,p-1\}, 
 \ k\in \mathbb N_0 \right\}
 $$
 is called the ring of $p$-adic integers in $\mathbb{Q}_{p}$. We know that $\mathbb{Z}_{p}$ is a compact set and $\mathbb{N}$ is
dense in $\mathbb{Z}_{p}$ (\cite{go}).

Throughout this article we shall consider all vector spaces and algebras taken over the
field $\mathbb{Q}_p$ (unless stated otherwise).

\section{Main results}

We are ready to present our results. For the rest of this work, $X$ will denote a nonempty subset of $\mathbb K$ without isolated points.
Let us fix two notations:
	$$
	C^1(X,\mathbb K)=\{f\colon X \to \mathbb K \colon f \text{ has continuous (classical) derivative on }X\},
	$$
	$$
	S^r(X,\mathbb K)=\{f\colon X \to \mathbb K \colon f \text{ is strictly differentiable of order }r \text{ on }X\}.
	$$

Our first result shows that, unlike the classical case, strict  differentiability is a stronger condition than simply having continuous derivative. 
An analogue to part (ii) of the result  for the classical case can be found in \cite{bbf}.

\begin{theorem}
\label{34} 
~\begin{itemize}
\item[(i)] The set $C^1(\qp,\qp)\setminus S^{1}(\qp,\qp)$ is $\mathfrak{c}$-lineable.
\item[(ii)] The set $S^{1}(\qp,\qp)\setminus S^{2}(\qp,\qp)$ is $\mathfrak{c}$-lineable. In general, $S^{r}(\qp,\qp)\setminus S^{r+1}(\qp,\qp)$ is $\mathfrak{c}$-lineable for any $r\geq 1$.
\end{itemize}
\end{theorem}

\begin{proof}	
\textit{(i).} Notice that $B_{p^{-2n}}(p^n)\subset S_{p^{-n}}(0)$ for every $n\in\mathbb N$, therefore $B_n\cap B_m=\emptyset$ if, and only if, $n\neq m$.
Also, let us define $f_N\colon \qp\rightarrow \qp$ for every $N\in\mathcal N$ as follows:
$$
f_N(x)=\begin{cases}
p^{2n} & \text{if }x\in B_{p^{-2n}}(p^n) \text{ with } n\in N,\\
0 & \text{otherwise}.
\end{cases}
$$
	First, notice that $f_N$ is locally constant outside $0$ and, thus, $f_N$ is differentiable everywhere except maybe at  $0$ with $f_N'(x)=0$ for every $x\in\qp\setminus\{0\}$.
	Moreover, we have 
	$$
	\left|\frac{f_N(x)-f_N(0)}{x}\right|_p=\left|\frac{f_N(x)}{x}\right|_p=\begin{cases}
	p^{-n} & \text{if } x\in B_{p^{-2n}}(p^n) \text{ with } n\in N,\\
	0 & \text{otherwise},
	\end{cases}
	$$
	i.e., $\left|\frac{f_N(x)-f_N(0)}{x}\right|_p\rightarrow 0$ as $x\rightarrow 0$.
	Hence, $f_N^\prime(0)=0$.
	Therefore, $f_N^\prime$ exists everywhere and is continuous since $f_N^\prime\equiv 0$, that is, $f_N\in C^1(\qp,\qp)$.
	
	It is enough to show that the family of functions $V=\{f_N\colon N\in\mathcal N \}$ is linearly independent over $\qp$ and the vector space generated by $V$, denoted by $\text{span}(V)$, satisfies $\text{span}(V)\setminus \{0\}\subset C^1(\qp,\qp)\setminus S^{1}(\qp,\qp)$.
	Take $f=\sum_{i=1}^m \alpha_i f_{N_i}$, where $\alpha_1,\ldots, \alpha_m\in\qp$, $N_1,\ldots, N_m\in\mathcal N$ are distinct and $m\in \mathbb N$.
	Assume that $f\equiv 0$ then, by taking $x=p^n$ with $n\in N_1^1\cap N_2^0\cap\cdots \cap N_m^0$, we have that $0=f(x)=\alpha_1 f_{N_1}(x)=\alpha_1 p^{2n}$, i.e., $\alpha_1=0$.
	Therefore, applying similar arguments, we arrive at $\alpha_i=0$ for every $i\in\{1,\ldots, m \}$.
	Assume now that $\alpha_i\neq 0$ for every $i\in\{1,\ldots, m \}$.
	Since $C^1(\qp,\qp)$ forms a vector space, we have that $f\in C^1(\qp,\qp)$.
	Moreover, by construction we have $f'\equiv 0$.
	It remains to prove that $f\notin S^1(\qp,\qp)$.
	To do so, take the sequences $$(x_n)_{n\in N_1^1\cap N_2^0\cap \cdots\cap N_m^0}=(p^n)_{n\in N_1^1\cap N_2^0\cap \cdots\cap N_m^0}$$ and $$(y_n)_{n\in N_1^1\cap N_2^0\cap \cdots\cap N_m^0}=(p^n-p^{2n})_{n\in N_1^1\cap N_2^0\cap \cdots\cap N_m^0}.$$
	Notice that both sequences converge to $0$, and $f(x_n)=\alpha_1 p^{2n}$ and $f(y_n)=0$ for each $n\in N_1^1\cap N_2^0\cap \cdots\cap N_m^0$.
	Hence, 
		$$
		\frac{f(x_n)-f(y_n)}{x_n-y_n}=\alpha_1,
		$$
	for every $n\in N_1^1\cap N_2^0\cap \cdots\cap N_m^0$.
	Since $\alpha_1\neq 0$, we have the desired result.
	
	\textit{(ii).} We will prove the case when $r=1$ since the general case can be easily deduced.
	For every nonempty subset $N$ of $\mathbb N$, let us define $g_N \colon \qp\rightarrow\qp$ as follows: for every $x=\sum_{n=m}^\infty a_n p^n\in\qp$, take
		$$
		g_N(x)=\sum_{n=0}^\infty b_n p^{2n},
		$$
	where
		$$
		b_n=\begin{cases}
		a_n & \text{if } n\in N,\\
		0 & \text{otherwise},
		\end{cases}
		$$
	for every $n\in\mathbb N_0$.
	For any $x,y\in \qp$, notice that 
		$$
		|g_N(x)-g_N(y)|_p\leq |x-y|_p^2.
		$$
	Hence, 
		$$
		\left|\frac{g_N(x)-g_N(z)}{x-z} \right|_p\leq |x-z|_p\rightarrow 0
		$$
	as $x\rightarrow z$ for any $z\in\qp$, that is, $g_N$ is differentiable and $g_N^\prime\equiv 0$.
	Moreover,
		$$
		\left|\frac{g_N(x)-g_N(y)}{x-y} \right|_p\leq |x-y|_p\rightarrow 0,
		$$
	when $(x,y)\rightarrow (z,z)$ for any $z\in\qp$, and where $(x,y)\in \nabla^2 \qp$.
	Thus, $g_N\in S^{1}(\qp,\qp)$.

	We will prove that the family of functions $W=\{g_N \colon N\in\mathcal N \}$ is linearly independent over $\qp$ and $\text{span}(W)\setminus\{0\}\subset S^{1}(\qp,\qp)\setminus S^{2}(\qp,\qp)$.
	It is easy to see that any linear combination of the functions in $W$ belongs to $S^{1}(\qp,\qp)$.
	Take now $g=\sum_{i=1}^k \beta_i g_{N_i}$, where $\beta_1,\ldots,\beta_k\in \qp$, $N_1,\ldots,N_k\in\mathcal N$ are distinct and $k\in\mathbb N$.
	Assume that $g\equiv 0$ then, by taking $x=p^{n}$, with $n\in N_1^1\cap N_2^0\cap \cdots\cap  N_k^0$ fixed, we have that $0=g(x)=\beta_1 p^{2n}$, i.e., $\beta_1=0$.
	By applying similar arguments we see that $\beta_i=0$ for every $i\in \{1,\ldots,k \}$.
	Therefore, assume that $\beta_i\neq 0$ for every $i\in \{1,\ldots,k \}$.
	For every $n\in  N_1^1\cap N_2^0\cap \cdots \cap N_k^0$, denote
		$$
		n_{+}=\min\{l\in N_1^1\cap N_2^0\cap \cdots \cap N_k^0 \colon l>n \}.
		$$
	Now consider the sequences $\overline x=(x_n)_{N_1^1\cap N_2^0\cap \cdots \cap  N_k^0}$, $\overline y=(y_n)_{N_1^1\cap N_2^0\cap \cdots \cap N_k^0}$ and $\overline z=(z_n)_{n\in N_1^1\cap N_2^0\cap \cdots \cap N_k^0}$ defined as $x_n=p^n$, $y_n=0$ and $z_n=p^{n}+p^{n_+}$ for every $n\in N_1^1\cap N_2^0\cap \cdots \cap N_k^0$.
	(Notice that the sequences $\overline{x}$, $\overline{y}$ and $\overline{z}$ converge to $0$.)
	Then, 
		\begin{align*}
			& \left|(y_n-z_n)^{-1}\right|_p \left|\frac{g(x_n)-g(y_n)}{x_n-y_n}-\frac{g(x_n)-g(z_n)}{x_n-z_n} \right|_p\\
			& = \left|(p^n+p^{n_+})^{-1}\right|_p \left|\frac{\beta_1 p^{2n}}{p^n}-\frac{\beta_1 p^{2n}-\beta_1 p^{2n}-\beta_1 p^{2n_{+}} }{p^n-p^n-p^{n_+}} \right|_p\\
			& = |\beta_1|_p \left|\frac{p^n-p^{n_+}}{p^n+p^{n_+}} \right|_p=|\beta_1|_p\neq 0,
		\end{align*}
	for every $n\in N_1^1\cap N_2^0\cap \cdots \cap  N_k^0$.
	However, $g^{\prime\prime}\equiv 0$.
	This finishes the proof.
\end{proof}

Let $\mathbb K$ be a non-Archimedean field with absolute value $|\cdot|$ that contains $\qp$.
For any $\alpha>0$, the space of Lipschitz functions from $X$ to $\mathbb K$ of order $\alpha$ is defined as
	$$
	\text{Lip}_\alpha (X,\mathbb K)=\{f\colon X\rightarrow \mathbb K \colon \exists M>0(|f(x)-f(y)|\leq M|x-y|^\alpha),\forall x,y\in X \}.
	$$
Let
	$$
	N^1(X,\mathbb K)=\{f\in S^1(X,\mathbb K)\colon f' \equiv 0 \}.
	$$
In view of \cite{sc}*{Exercise~63.C} we have
	$$
	N^1(\zp,\mathbb K)\setminus \bigcup_{\alpha>1} \text{Lip}_\alpha (\zp, \mathbb K)\neq \emptyset.
	$$
	
To prove the next theorem, we need a  characterization of the spaces $N^1(\zp,\mathbb K)$ and $\text{Lip}_\alpha (\zp,\mathbb K)$ from \cite{sc}. 
For this we will denote by $(e_n)_{n\in \mathbb N_0}$ the van der Put base of $C(\zp,\mathbb K)$, which is given by $e_0\equiv 1$ and $e_n$ is  the characteristic function of $\{x\in\zp \colon |x-n|_p<1/n \}$ for every $n\in \mathbb N$.

\begin{proposition}\label{thm:1}~\begin{itemize}
	\item[(i)] 	Let $f=\sum_{n=0}^\infty a_n e_n \in C(\zp,\mathbb K)$.
	Then $f\in N^1(\zp,\mathbb K)$ if and only if  $(|a_n| n)_{n\in \mathbb N_0}$ converges to $0$ (see \cite{sc}*{Theorem~63.3}).\label{thm:63.3}
	
	\item[(ii)]
	Let $f=\sum_{n=0}^\infty a_n e_n \in C(\zp,\mathbb K)$ and $\alpha>0$.
	Then $f\in \text{Lip}_\alpha (\zp,\mathbb K)$ if and only if 
	$$
	\sup \{|a_n| n^\alpha \colon n\in\mathbb N_0 \}<\infty
	$$ 
	(see \cite{sc}*{Exercise~63.B}).\label{ex:63.B}
\end{itemize}
\end{proposition} 

The next result shows that there is a vector space of maximum dimension of  strictly differentiable functions with zero derivative that are not Lipschitzian.
Our next three results can be compared with some results obtained in \cites{jms, bbfg,bq} for the classical case.

\begin{theorem}
	The set $N^1(\mathbb{Z}_{p},\mathbb K)\setminus\bigcup_{\alpha>1}{\rm{Lip}}_{\alpha }(\mathbb{Z}_{p},\mathbb K)$ is $\mathfrak{c}$-lineable (as a $\mathbb K$-vector space).
\end{theorem}

\begin{proof}
	Fix $n_1\in\mathbb N$ and take $B_1=\{x\in \mathbb Z_p\colon |x-n_1|_p<1/n_1 \}$.
	Since $\mathbb Z_p$ and $B_1$ are clopen sets, we have that $\mathbb Z_p\setminus B_1$ is open and also nonempty.
	Therefore there exists an open ball $D_1\subset \mathbb Z_p\setminus B_1$.
	Furthermore, as $\mathbb N$ is dense in $\mathbb Z_p$, there exists $n_2\in \mathbb N\setminus \{n_1\}$ such that $B_2=\{ x\in \mathbb Z_p\colon |x-n_2|_p<1/n_2 \}\subset D_1$.
	By recurrence, we can construct a  set  $M=\{n_k\colon k\in\mathbb N \}\subset \mathbb N$ such that the balls $B_k=\{x\in\mathbb Z_p\colon |x-n_k|_p<1/n_k \}$ are pairwise disjoint.
	
	Let $\sigma\colon \mathbb N_0\rightarrow M$ be the increasing bijective function and let $(m_n)_{n\in\mathbb N}\subset \mathbb N$ be an increasing sequence such that $p^{-m_n}n\rightarrow 0$ and $p^{-m_n}n^\alpha\rightarrow\infty$ for every $\alpha>1$ when $n\rightarrow\infty$. (This can be done for instance by taking $m_n=\left\lfloor \ln(n \ln(n))/\ln(p) \right\rfloor$.)
	For every $N\in\mathcal N_0$, define $f_N\colon \mathbb Z_p\rightarrow \mathbb K$ as
		$$
		f_N=\sum_{n=0}^\infty 1_N(n) p^{m_{\sigma(n)}} e_{\sigma(n)}.
		$$
	Since every $N\in\mathcal N_0$ is infinite, we have that $|1_N(n)| p^{-m_{\sigma(n)}} \sigma(n)\rightarrow 0$ when $n\rightarrow \infty$ and 
		$$
		\{|1_N(n)| p^{-m_{\sigma(n)}} \sigma(n)^\alpha \colon n\in\mathbb N_0 \}
		$$ 
	is unbounded for every $\alpha>1$.
	Hence, by Theorem~\ref{thm:1}, we have $f_N\in N^1(\mathbb{Z}_{p},\mathbb K)\setminus\bigcup_{\alpha>1}{\rm{Lip}}_{\alpha }(\mathbb{Z}_{p},\mathbb K)$ for every $N\in\mathcal N_0$.
	
	We will prove now that the functions in $V=\{f_N\colon N\in\mathcal N_0 \}$ are linearly independent over $\mathbb K$ and such that any nonzero linear combination of $V$ over $\mathbb K$ is contained in $N^1(\mathbb{Z}_{p},\mathbb K)\setminus\bigcup_{\alpha>1}{\rm{Lip}}_{\alpha }(\mathbb{Z}_{p},\mathbb K)$.
	Take $f=\sum_{i=1}^r a_i f_{N_i}$, where $a_1,\ldots, a_r\in\mathbb K$, $N_1,\ldots, N_r\in\mathcal N_0$ are distinct and $r\in\mathbb N$.
	Assume that $f\equiv 0$. Fix $n\in N_1^1\cap N_2^0\cap \cdots\cap N_r^0$ and take $x\in\mathbb Z_p$ such that $x\in B_{\sigma(n)}$, then $0=f(x)=a_1 p^{m_{\sigma(n)}}$, i.e., $a_1=0$.
	By applying similar arguments we have $a_i=0$ for every $i\in\{1,\ldots,r \}$.
	Assume for the rest of the proof that $a_i\neq 0$ for every $i\in\{1,\ldots,r \}$.
	Notice that $f=\sum_{n=0}^\infty \left(\sum_{i=1}^r a_i 1_{N_i}\right)(n) p^{m_{\sigma(n)}} e_{\sigma(n)}$, where $\left|\left(\sum_{i=1}^r a_i 1_{N_i}\right)(n) p^{m_{\sigma(n)}}\right|\leq |p^{m_{\sigma(n)}}| \max\{|a_i|\colon i\in\{1,\ldots,r \} \}=p^{-m_{\sigma(n)}} \max\{|a_i|\colon i\in\{1,\ldots,r \} \}$.
	Therefore $\left|\left(\sum_{i=1}^r a_i 1_{N_i}\right)(n) p^{m_{\sigma(n)}}\right| \sigma(n)\rightarrow 0$ when $n\rightarrow \infty$.
	Finally, as $N_1^1\cap N_2^0\cap \cdots\cap N_r^0$ is infinite, we have that 
		\begin{align*}
			& \left\{\left|\left(\sum_{i=1}^r a_i 1_{N_i}\right)(n) p^{m_{\sigma(n)}}\right| \sigma(n)^\alpha \colon n\in N_1^1\cap N_2^0\cap \cdots\cap N_r^0 \right\}\\
			& = \{|a_1| p^{-m_{\sigma(n)}} \sigma(n)^\alpha \colon n\in N_1^1\cap N_2^0\cap \cdots\cap N_r^0 \}
		\end{align*}
	is unbounded for every $\alpha>1$.
\end{proof}

The next lineability result can be considered as a non-Archimedean counterpart of  \cite{g7}*{Theorem 6.1}.
To prove it we will make use of the following lemma.
(For more information on the usage of the lemma see \cite[Section~5]{fmrs}.)
In order to understand it, let us consider the functions $x\mapsto (1+x)^\alpha$ where $x\in p\mathbb Z_p$ and $\alpha\in \mathbb Z_p$.
It is well known that $(1+x)^\alpha$ is defined analytically by $(1+x)^\alpha = \sum_{i=0}^\infty {\alpha \choose i} x^i$.
Moreover $x\mapsto (1+x)^\alpha$ is well defined (see \cite[Theorem~47.8]{sc}), differentiable with derivative $\alpha(1+x)^{\alpha-1}$, and $x\mapsto (1+x)^\alpha$ takes values in $\mathbb Z_p$ (in particular $(1+x)^\alpha=1+y$ for some $y\in p\mathbb Z_p$, see \cite[Theorem~47.10]{sc}).

\begin{lemma}\label{lem:1}
	If $\alpha_1,\ldots,\alpha_n\in \mathbb Z_p\setminus \{0\}$ are distinct, with $n\in \mathbb N$, then every linear combination $\sum_{i=1}^n \gamma_i (1+x)^{\alpha_i}$, with $\gamma_i\in \mathbb Q_p\setminus \{0\}$ for every $1\leq i\leq n$, is not constant.
\end{lemma}

\begin{theorem}
	\label{16} The set of everywhere differentiable functions $\mathbb{Q}_{p}\rightarrow \mathbb{Q}_{p}$
	which derivative is unbounded is strongly $\mathfrak{c}$-algebrable.
\end{theorem}

\begin{proof}
	Let $\mathcal H$ be a Hamel basis of $\mathbb Q_p$ over $\mathbb Q$ contained in $p\mathbb Z_p$, and for each $\beta \in \mathbb Z_p\setminus\{0\}$ define $f_\beta : \mathbb Q_p\to \mathbb Q_p$ by
		$$
		f_\beta (x)=\begin{cases}
		p^{-n} (1+y)^\beta & \text{if } x=\sum_{k=-n}^0 a_k p^k +y, \text{ where } a_n\neq 0,\ n\in \mathbb N_0 \text{ and } y\in p\zp,\\
		0 & \text{otherwise}.
		\end{cases}
		$$
	The function $f_\beta$ is differentiable everywhere for any $\beta\in \mathbb Z_p$.
	Indeed, firstly we have that $f_\beta$ is locally constant on $p\mathbb Z_p$ as $f_\beta\restriction p\mathbb Z_p \equiv 0$.
	Lastly it remains to prove that $f_\beta$ is differentiable at $x_0=\sum_{k=-n}^0 a_k p^k +y_0$, i.e., the limit
		\begin{equation}\label{equ:2}
			\lim_{x\rightarrow x_0} \frac{p^{-m}(1+y)^\beta-p^{-n}(1+y_0)^\beta}{x-x_0}
		\end{equation}
	exists, where the values $x$ are of the form $x=\sum_{k=-m}^0 b_k p^k +y$.
	Moreover, as $x$ approaches $x_0$ in the limit \eqref{equ:2}, we can assume that $x=\sum_{k=-n}^0 a_k p^k +y$.
	Hence, the limit in \eqref{equ:2} can be simplified to
		\begin{align*}
			\lim_{x\rightarrow x_0} \frac{p^{-m}(1+y)^\beta-p^{-n}(1+y_0)^\beta}{x-x_0} & =\lim_{x\rightarrow x_0} \frac{p^{-n}(1+y)^\beta-p^{-n}(1+y_0)^\beta}{x-x_0}\\
			& =p^{-n} \lim_{y\rightarrow y_0} \frac{(1+y)^\beta-(1+y_0)^\beta}{y-y_0}=p^{-n}\beta(1+y_0)^{\beta-1}.
		\end{align*}
	In particular the derivative of $f_\beta$ is given by
		$$
		f_\beta^\prime (x)=\begin{cases}
		p^{-n} \beta (1+y)^{\beta-1} & \text{if } x=\sum_{k=-n}^0 a_k p^k +y, \text{ where } a_n\neq 0,\ n\in \mathbb N_0, y\in p\zp,\\
		0 & \text{otherwise},
		\end{cases}
		$$
	and it is unbounded since $\lim_{n\rightarrow \infty} |p^{-n}\beta(1+y)^{\beta-1}|_p=\lim_{n\rightarrow \infty} p^{n}|\beta(1+y)^{\beta-1}|_p=|\beta|_p\lim_{n\rightarrow \infty} p^{n}=\infty$, where we have used the fact that $\beta\neq 0$.
	
	Let $h_1,\ldots,h_m \in \mathcal H$ be distinct and take $P$ a polynomial in $m$ variables with coefficients in $\mathbb Q_p\setminus \{0\}$ and without free term, that is, $P(x_1,\ldots,x_m)=\sum_{r=1}^d \alpha_r x_1^{k_{r,1}}\cdots x_m^{k_{r,m}}$, where $d\in \mathbb N$, $\alpha_i\in\qp\setminus \{0\}$ for every $1\leq r\leq d$, $k_{r,i}\in\mathbb N_0$ for every $1\leq r\leq d$ and $1\leq i\leq m$ with $\overline{k}_r:=\sum_{i=1}^m k_{r,i}\geq 1$, and the $m$-tuples $(k_{r,1},\ldots,k_{r,m})$ are pairwise distinct.
	Assume also without loss of generality that $\overline{k}_1\geq \cdots\geq \overline{k}_d$.
	We will prove first that $P(f_{h_1},\ldots,f_{h_m})\not\equiv 0$, i.e., the functions in $\{f_h \colon h\in \mathcal  H \}$ are algebraically independent.
	Notice that $P(f_{h_1},\ldots,f_{h_m})$ is of the form
		$$
		\begin{cases}
		\sum_{r=1}^d p^{-n\overline{k}_r} \alpha_r (1+y)^{\beta_r} & \text{if } x=\sum_{k=-n}^0 a_k p^k +y, \text{ with } a_n\neq 0,\ n\in \mathbb N_0,  y\in p\zp,\\
		0 & \text{otherwise},
		\end{cases}
		$$
	where the exponents $\beta_r:=\sum_{i=1}^m k_{r,i} h_i$ belong to $p\zp\setminus\{0\}$ and are pairwise distinct because $\mathcal H$ is Hamel basis of $\mathbb Q_p$ over $\mathbb Q$ contained in $p\zp$, $k_{r,i}\in\mathbb N_0$, $\overline{k}_r\neq 0$ and the numbers $h_1,\ldots,h_m$ as well as the $m$-tuples $(k_{r,1},\ldots,k_{r,m})$ are pairwise distinct.	
	Fix $n\in\mathbb N_0$.
	Since $p^{-n\overline{k}_r} \alpha_r \neq 0$ for every $1\leq r\leq d$, by Lemma~\ref{equ:2}, there is $y\in p\zp$ such that $\sum_{r=1}^d p^{-n\overline{k}_r} \alpha_r (1+y)^{\beta_r}\neq 0$.
	Hence, by taking $x=p^{-n}+y$, we have $P(f_{h_1},\ldots,f_{h_m})(x)\neq 0$.
	
	Finally, let us prove that $P(f_{h_1},\ldots,f_{h_m})'$ exists and it is unbounded.
	Clearly $P(f_{h_1},\ldots,f_{h_m})$ is differentiable and the derivative is given by
		\begin{equation}\label{equ:6}
			\begin{cases}
			\displaystyle \sum_{r=1}^d p^{-n\overline{k}_r} \alpha_r \beta_r (1+y)^{\beta_r-1} & \displaystyle \text{if } x=\sum_{k=-n}^0 a_k p^k +y, \text{ with } a_n\neq 0,\ n\in \mathbb N_0, y\in p\zp,\\
			0 & \text{otherwise}.
			\end{cases}
		\end{equation}
	Notice that $\beta_r\neq 1$ for every $1\leq r\leq d$ since $\beta_r\in p\zp$.
	We will now rewrite formula \eqref{equ:6} in order to simplify the proof.
	Notice that if some of the exponents $\overline{k}_r$ are equal, i.e., for instance $\overline{k_i}=\cdots=\overline{k}_j$ for some $1\leq i<j\leq d$, then $p^{-n\overline{k}_i}$ is a common factor in each summand $p^{-n\overline{k}_i} \alpha_i \beta_i (1+y)^{\beta_i-1},\ldots,p^{-n\overline{k}_j} \alpha_j \beta_j (1+y)^{\beta_j-1}$.
	Therefore, $P(f_{h_1},\ldots,f_{h_m})^\prime(x)$ can also be written as
		\begin{equation}\label{equ:7}
			\begin{cases}
			\displaystyle \sum_{q=1}^{\widetilde{d}} p^{-n\widetilde{k}_q} \sum_{s=1}^{m_q}\alpha_{q,s}\beta_{q,s} (1+y)^{\beta_{q,s}-1} & \displaystyle \text{if } x=\sum_{k=-n}^0 a_k p^k +y, \text{ where } a_n\neq 0,\\
			& \displaystyle n\in \mathbb N_0 \text{ and } y\in p\zp,\\
			0 & \displaystyle \text{otherwise},
			\end{cases}
		\end{equation}
	where $\tilde{d}\in \mathbb{N}$, the $\widetilde{k}_q$'s represent the common exponents $p^{-n\overline{k}_i}$ with $\widetilde{k}_1>\cdots>\widetilde{k}_{\widetilde{d}}$, and $\alpha_{q,s}$ and $\beta_{q,s}$ are the corresponding terms $\alpha_r$ and $\beta_r$, respectively.
	Now, since $\alpha_{1,s} \beta_{1,s}\neq 0$ and the exponents $\beta_{1,s}-1$ are pairwise distinct and not equal to $0$ for every $1\leq s\leq m_1$, by Lemma~\ref{equ:2}, there exists $y_1\in p\zp$ such that $\sum_{s=1}^{m_1}\alpha_{1,s}\beta_{1,s} (1+y_1)^{\beta_{1,s}-1}\neq 0$.
	Take the sequence $(x_n)_{n=1}^\infty$ defined by $x_n=p^{-n}+y_1$.
	Since $\widetilde{k}_1>\cdots>\widetilde{k}_{\widetilde{d}}$, there exists $n_0\in \mathbb N$ such that
		\begin{align*}
			|P(f_{h_1},\ldots,f_{h_m})^\prime (x_n)|_p & =\left|\sum_{q=1}^{\widetilde{d}} p^{-n\widetilde{k}_q} \sum_{s=1}^{m_q}\alpha_{q,s}\beta_{q,s} (1+y_1)^{\beta_{q,s}-1} \right|_p\\
			& =\left|p^{-n\widetilde{k}_1}\sum_{s=1}^{m_1 }\alpha_{1,s} \beta_{1,s} (1+y_1)^{\beta_{1,s}-1} + \right.\\ 
			&\left. + \sum_{q=2}^{\widetilde{d}} p^{-n\widetilde{k}_r} \sum_{s=1}^{m_q}\alpha_{q,s} \beta_{q,s} (1+y_1)^{\beta_{q,s}-1} \right|_p\\
			& =\left|p^{-n\widetilde{k}_1}\sum_{s=1}^{m_1 }\alpha_{1,s} \beta_{1,s} (1+y_1)^{\beta_{1,s}-1}\right|_p\\
			& =p^{n\widetilde{k}_1}\left|\sum_{s=1}^{m_1 }\alpha_{1,s} \beta_{1,s} (1+y_1)^{\beta_{1,s}-1}\right|_p,
		\end{align*}
	for every $n\geq n_0$. 
	Hence, $\lim_{n\to \infty} |P(f_{h_1},\ldots,f_{h_m})^\prime (x_n)|_p=\infty$.
\end{proof}

The reader may have noticed that the functions in the proof of Theorem~\ref{16} have unbounded derivative but the derivative is bounded on each ball of $\qp$.
The following result (which proof is a modification of the one in Theorem~\ref{16}) shows that we can obtain a similar optimal result when the derivative is unbounded on each ball centered at a fixed point $a\in \qp$.
The functions will not be differentiable at $a$.

\begin{corollary}
\label{15} Let $a\in \qp$.
The set of continuous functions $\mathbb{Q}_{p}\rightarrow \mathbb{Q}_{p}$ that are differentiable except at $a$ and which derivative is unbounded on $\qp\setminus (a+\zp)$ and on $(a+\zp)\setminus\{a\}$ is strongly $\mathfrak{c}$-algebrable.
\end{corollary}

\begin{proof}
	Fix $a\in \qp$.
	Let $\mathcal H$ be a Hamel basis of $\mathbb Q_p$ over $\mathbb Q$ contained in $p\mathbb Z_p$.
	For every $\beta\in \zp\setminus\{0\}$ take the function $f_\beta$ defined in the proof of Theorem~\ref{16} and also define $g_\beta\colon \qp\rightarrow \qp$ by
		$$
		g_\beta(x)=\begin{cases}
		p^{n}[p^{-n^2}(x-a)]^\beta & \text{if } x\in \overline{B}_{p^{-(n^2+1)}}(a+p^{n^2}) \text{ for some } n\in \mathbb N,\\
		0 & \text{otherwise}.
		\end{cases}
		$$
	Notice that by applying the change of variable $y=x-a$ we can assume, without loss of generality, that $a=0$.
	Since for any $x\in \overline{B}_{p^{-(n^2+1)}}(p^{n^2})$ with $n\in\mathbb N$, $x$ is of the form $p^{n^2}+\sum_{k=n^2+1}^\infty a_k p^k$ with $a_{k}\in \{0,1,\ldots,p-1\}$ for every integer $k\geq n^2+1$, we have that $p^{-n^2}x=1+\sum_{k=n^2+1}^\infty a_k p^{k-n^2}\in 1+p\zp$.
	Thus $g_\beta$ is well defined.
	Now, for every $\beta\in \zp\setminus\{0\}$, let $F_{\beta}:=f_\beta+g_\beta$.
	It is easy to see that $F_\beta$ is differentiable at every $x\in \qp\setminus\{0\}$ and, in particular, continuous on $\qp\setminus\{0\}$.
	Let us prove now that $F_\beta$ is continuous at $0$.
	Fix $\varepsilon>0$ and take $n\in\mathbb N$ such that $p^{-n}<\varepsilon$. 
	If $|x|_p<p^{1-n^2}$, then 
		\begin{align*}
			|F_\beta(x)| & =\begin{cases}
			|p^n (p^{-n^2}x)^\beta|_p & \text{if } x\in \overline{B}_{p^{-(n^2+1)}}(p^{n^2}),\\
			0 & \text{otherwise},
			\end{cases}\\
			& =\begin{cases}
			p^{-n} & \text{if } x\in \overline{B}_{p^{-(n^2+1)}}(p^{n^2}),\\
			0 & \text{otherwise}.
			\end{cases}
		\end{align*}
	In any case, $|F_\beta(x)|<\varepsilon$.
	Hence $F_\beta$ is continuous.
	Moreover, $F_\beta$ is not differentiable at $0$.
	Indeed, by considering the sequence $(x_n)_{n=1}^\infty=(p^{n^2}+p^{n^2+1})_{n=1}^\infty$ which converges to $0$ when $n\to \infty$ we have
		\begin{align*}
			\lim_{n\to \infty} \frac{|F_\beta(x_n)-F_\beta(0)|_p}{|x_n|_p} & =\lim_{n\to \infty} \frac{|p^{n} [p^{-n^2}(p^{n^2}+p^{n^2+1})]^\beta|_p}{|p^{n^2}+p^{n^2+1}|_p}\\
			& =\lim_{n\to \infty} \frac{p^{-n} |(1+p)^\beta|_p}{p^{-n^2}}=\lim_{n\to \infty} p^{n^2-n}=\infty.
		\end{align*}
	In particular, by the chain rule, the derivative of $F_\beta$ on $\qp\setminus \{0\}$ is as follows
		$$
		\hspace*{-1cm} F_\beta^\prime(x)=\begin{cases}
		p^{-n} \beta (1+y)^{\beta-1} & \displaystyle \text{if } x=\sum_{k=n}^0 a_k p^k +y, \text{ with } a_n\neq 0,\ n\in \mathbb Z\setminus \mathbb N,  y\in p\zp,\\
		p^{n-n^2}\beta(p^{-n^2}x)^{\beta-1} & \text{if } x\in \overline{B}_{p^{-(n^2+1)}}(p^{n^2}) \text{ for some } n\in \mathbb N,\\
		0 & \text{otherwise}.
		\end{cases}
		$$
		
	By the proof of Theorem~\ref{16} the functions in the set $V=\{F_h \colon h\in \mathcal H \}$ are algebraically independent, and every function in the algebra $\mathcal A$ generated by $V$ over $\qp$ that is not the $0$ function is continuous, differentiable on $\qp\setminus \{0\}$ and has unbounded derivative on $\qp\setminus \zp$.
	It remains to prove that any nonzero algebraic combination in $V$ has unbounded derivative on $\zp\setminus\{0\}$.
	To do so, let $h_1,\ldots,h_m \in \mathcal H$ be distinct and take $P$ a polynomial in $m$ variables with coefficients in $\mathbb Q_p\setminus \{0\}$ and without free term.
	Then, by the chain rule, $P(f_{h_1},\ldots,f_{h_m})^\prime$ on $\overline{B}_{p^{-(n^2+1)}}(p^{n^2})$ is of the form
		$$
		p^{-n^2}\sum_{q=1}^d p^{nk_q} \sum_{s=1}^{m_q} \alpha_{q,s} \beta_{q,s} (p^{-n^2}x)^{\beta_{q,s}-1},
		$$
	see the proof of Theorem~\ref{16}.
	Assume without loss of generality that $k_1<\cdots<k_d$.
	Since $\alpha_{1,s} \beta_{1,s}\neq 0$ and the exponents $\beta_{1,s}- 1$ are pairwise distinct and not $0$ for every $1\leq s\leq m_1$, by Lemma~\ref{lem:1}, there exists $y_1\in p\zp$ such that $\sum_{s=1}^{m_1} \alpha_{1,s} \beta_{1,s}(1+y_1)^{\beta_{1,s}-1}\neq 0$.
	For every $n\in\mathbb N$, take $x_n=p^{n^2}(1+y_1)$.
	Then, notice that there exists $n_0\in\mathbb N$ such that 
		\begin{align*}
			|P(F_{h_1},\ldots,F_{h_m})^\prime (x_n)|_p & =\left|p^{-n^2}\sum_{q=1}^{d} p^{nk_q} \sum_{s=1}^{m_q}\alpha_{q,s}\beta_{q,s} (p^{-n^2}x_n)^{\beta_{q,s}-1} \right|_p\\
			& =p^{n^2}\left|\sum_{q=1}^{d} p^{nk_q} \sum_{s=1}^{m_q}\alpha_{q,s}\beta_{q,s} (1+y_1)^{\beta_{q,s}-1} \right|_p\\
			& =p^{n^2}\left|p^{nk_1}\sum_{s=1}^{m_1 }\alpha_{1,s} \beta_{1,s} (1+y_1)^{\beta_{1,s}-1} + \right. \\
			& \left. + \sum_{q=2}^{d} p^{nk_r} \sum_{s=1}^{m_q} \alpha_{q,s} \beta_{q,s} (1+y_1)^{\beta_{q,s}-1} \right|_p\\
			& =p^{n^2}\left|p^{nk_1}\sum_{s=1}^{m_1 }\alpha_{1,s} \beta_{1,s} (1+y_1)^{\beta_{1,s}-1}\right|_p\\
			& =p^{n^2-nk_1}\left|\sum_{s=1}^{m_1 }\alpha_{1,s} \beta_{1,s} (1+y_1)^{\beta_{1,s}-1}\right|_p,
			\end{align*}
	for every $n\geq n_0$. 
	Therefore $\lim_{n\to \infty} |P(F_{h_1},\ldots,F_{h_m})^\prime (x_n)|_p=\infty$.
\end{proof}

In Corollary~\ref{15} we can replace unbounded derivative with being not Lipschitzian although the conclusion is weaker in terms of lineability as it is shown in the following proposition.

\begin{proposition}
	\label{26} The set of continuous functions $\mathbb{Q}_{p}\rightarrow \mathbb{Q}_{p}$ which are differentiable except at $0$, with bounded derivative on $\qp\setminus\{0\}$ and not Lipschitzian of order $\alpha>0$ is $\mathfrak{c}$-lineable and 1-algebrable.
\end{proposition}

\begin{proof}
	Let us prove first the lineability part.
	For any $N\in\mathcal N$, let $f_N\colon \qp\rightarrow \qp$ be:
	$$
	f_N(x)=\begin{cases}
		p^{n} & \text{if } x\in S_{p^{-n^2}}(0) \text{ and } n\in N,\\
		0 & \text{otherwise}.
	\end{cases}
	$$
	For every $x\in\qp\setminus\{0\}$, it is clear that there exists a neighborhood of $x$ such that $f_N$ is constant since the spheres are open sets.
	Thus, $f_N$ is locally constant on $\qp\setminus\{0\}$ which implies that $f_N$ is continuous, differentiable on $\qp\setminus\{0\}$ and $f_N^\prime(x)=0$ for every $x\in \qp\setminus\{0\}$.
	Moreover, it is easy to see that $f_N$ is continuous at $0$.
	However, $f_N$ is not differentiable at $0$.
	Indeed, take $x_n=p^{n^2}$ for every $n\in N$.
	It is clear that the sequence $(x_n)_{n\in N}$ converges to $0$ and also, for every $\alpha>0$,
	$$
	\frac{|f_N(x_n)|_p}{|x_n|_p^\alpha}=p^{\left(-1+\alpha n\right) n}\rightarrow \infty,
	$$
	when $n\in N$ tends to infinity.
	Therefore $f_N$ is not differentiable at $0$.
	Furthermore, notice that for any $M>0$ there are infinitely many $x\in \zp$ such that $|f_N(x)|_p> M|x|_p^\alpha$.
	Hence $f_N$ is not Lipschitzian of order $\alpha>0$.
	
	It remains to prove that the functions in $V=\{f_N\colon N\in\mathcal N \}$ are linearly independent over $\qp$ and such that the functions in $\text{span}(V)\setminus\{0\}$ are differentiable except at $0$, with bounded derivative on $\qp\setminus\{0\}$ and not Lipschitzian of order $\alpha>0$.
	Let $f=\sum_{i=1}^m \alpha_i f_{N_i}$, where $\alpha_1,\ldots, \alpha_m\in \qp$, $N_1,\ldots, N_m\in \mathcal N$ are distinct and $m\in \mathbb N$.
	Assume that $f\equiv 0$ and take $n\in N_1^1\cap N_2^0\cap \cdots \cap N_m^0$.
	Then $0=f(p^{n^2})=\alpha_1 p^{n}$ which implies that $\alpha_1=0$.
	Applying similar arguments we have that $\alpha_i=0$ for every $i\in \{1,\ldots, m \}$.
	Finally, assume that $\alpha_i\neq 0$ for every $i\in \{1,\ldots, m \}$.
	It is clear that $f$ is continuous on $\qp$ and differentiable on $\qp\setminus\{0\}$ with $f_N^\prime(x)=0$ for every $x\in \qp\setminus \{0\}$.
	Let $x_n=p^{n^2}$ for every $n\in N_1^1\cap N_2^0\cap \cdots \cap N_m^0$ and notice that $(x_n)_{n\in N_1^1\cap N_2^0\cap \cdots \cap N_m^0}$ converges to $0$.
	Moreover, for every $\alpha>0$,
	$$
	\frac{|f(x_n)|_p}{|x_n|_p^\alpha}=|\alpha_1|_p p^{\left(-1+\alpha n\right) n}\rightarrow \infty,
	$$
	when $n\in N_1^1\cap N_2^0\cap \cdots \cap N_m^0$ tends to infinity.
	Hence, $f$ is not differentiable at $0$, and also for every $M>0$ there are infinitely many $x\in \zp$ such that $|f(x)|_p> M|x|_p^\alpha$.
	
	For the algebrability part, let $g\colon \qp\rightarrow\qp$ be defined as:
	$$
	g(x)=\begin{cases}
		p^n & \text{if } x\in S_{p^{-n^2}}(0) \text{ for some } n\in\mathbb N,\\
		0 & \text{otherwise}.
	\end{cases}
	$$
	By applying similar arguments used in the first part of the proof, we have that $g$ is continuous, differentiable on $\qp\setminus\{0\}$ with $g^\prime(x)=0$ for every $x\in \qp\setminus\{0\}$ and not Lipschitzian of order $\alpha>0$.
	To finish the proof, let $G=\beta g^k$ where $\beta\in\qp\setminus\{0\}$ and $k\in\mathbb N$.
	It is obvious that $G$ is continuous, differentiable on $\qp\setminus\{0\}$ and $G^\prime(x)=0$ for every $x\in\qp\setminus\{0\}$.
	Now, let $x_n=p^{n^2}$ for every $n\in\mathbb N$.
	It is easy to see that $(x_n)_{n\in\mathbb N}$ converges to $0$ and
	$$
	\frac{|G(x_n)|_p}{|x_n|_p^\alpha}=|\beta|_p p^{\left(-k+\alpha n\right) n}\rightarrow \infty,
	$$
	when $n\rightarrow\infty$.
\end{proof}

Let $\mathcal B$ be the $\sigma$-algebra of all Borel subsets of $\zp$ and $\mu$ be any non-negative real-valued Haar measure on the measurable space $(\zp,\mathcal B)$.
In particular, if $\mu$ is normalized, then $\mu\left(x+p^n\zp\right)=p^{-n}$ for any $x\in \zp$ and $n\in\mathbb N$.
For the rest of the paper $\mu$ will denote a non-negative real-valued Haar measure on $(\zp,\mathcal B)$.
As usual, a Borel subset $B$ of $\zp$ is called a null set for $\mu$ provided that $\mu(B)=0$.
We also say that a Borel subset of $\zp$ is a full set for $\mu$ if $\zp\setminus B$ is a null set. (See \cite[Section~2.2]{Fo} for more details on the Haar measure.)

It is easy to see that the singletons of $\zp$ are null sets for any Haar measure $\mu$ on $(\zp,\mathcal B)$.
Therefore Proposition~\ref{26} states in particular that, for any Haar measure $\mu$ on $(\zp,\mathcal B)$, the set of continuous functions $\zp\rightarrow\qp$ that are differentiable except on a null set for $\mu$ of cardinality $1$, with bounded derivative elsewhere, is $\mathfrak{c}$-lineable.
The following result shows that a similar version can be obtained when we consider null sets of cardinality $\mathfrak{c}$ for any Haar measure $\mu$ on $(\zp,\mathcal B)$.
In order to prove it, we recall the following definition and result from probability theory.

\begin{definition}\label{def:1}
	Let $(\Omega,\mathcal F,P)$ be a probability space and $Y$ be a measurable real-valued function on $\Omega$. We say that $Y$ is a random variable.
\end{definition}

\begin{theorem}\normalfont{(Strong law of large numbers, see \cite[Theorem~22.1]{Bi}).}
	Let $(Y_n)_{n\in \mathbb N_0}$ be a sequence of independent and identically distributed real-valued random variables on a probability space $(\Omega,\mathcal F,P)$ such that, for each $n\in\mathbb N_0$, $E[Y_n]=m$ for some $m\in\mathbb R$ (where $E$ denotes the expected value).
	Then
	$$
	P\left(\left\{x\in \Omega\colon \exists\lim_{n\rightarrow \infty} \frac{\sum_{k=0}^{n-1} Y_k(x)}{n}=m \right\} \right)=1.
	$$
\end{theorem}

\begin{theorem}\label{thm:2}
	Let $\mu$ be a Haar measure on $(\zp,\mathcal B)$.
	The set of continuous functions $\zp\rightarrow\qp$ which are differentiable on a full set for $\mu$ with bounded derivative but not differentiable on the complement having cardinality $\mathfrak c$ is $\mathfrak c$-lineable.
\end{theorem}

\begin{proof}
	We will prove the result for $\mu$ the normalized Haar measure on $(\zp,\mathcal B)$ since any null set for the normalized Haar measure is also a null set for any other non-negative real-valued Haar measure on $(\zp,\mathcal B)$.
	This is an immediate consequence of Haar's Theorem which states that Haar measures are unique up to a positive multiplicative constant (see \cite[Theorem~2.20]{Fo}).

	Let $f\colon \zp\rightarrow \zp$ be defined as follows: for every $x=\sum_{i=0}^\infty x_i p^i\in \zp$, we have
		\begin{equation}\label{equ:3}
			f(x)=\begin{cases}
				x & \text{if } (x_{2i},x_{2i+1})\neq (0,0) \text{ for all } i\in \mathbb N_0,\\
				\displaystyle \sum_{i=0}^{2n+1} x_i p^i & \text{if } (x_{2i},x_{2i+1})\neq (0,0) \text{ for all } i\leq n \\
				 & \text{with } n\in \mathbb N_0 \text{ and } (x_{2n+2},x_{2n+3})= (0,0),\\
				0 & \text{if } (x_0,x_1)=(0,0).
			\end{cases}
		\end{equation}
	The function $f$ is continuous.
	Indeed, let $x=\sum_{i=0}^\infty x_i p^i\in \zp$ and fix $\varepsilon>0$.
	Take any $m\in \mathbb N_0$ such that $p^{-(2m+1)}<\varepsilon$.
	Then for any $y\in \zp$ such that $|x-y|_p<p^{-(2m+1)}$ we have that $y$ is of the form $\sum_{i=0}^{2m+1} x_i p^i + \sum_{i=2m+2}^\infty y_i p^i$.
	Hence, notice that in any possible case of $x$ given in \eqref{equ:3}, we have that $|f(x)-f(y)|_p< p^{-(2m+1)}<\varepsilon$.
	
	Let us define, for every $i\in \mathbb N_0$, the random variables $Y_i\colon \zp\rightarrow \{0,1\}$ in the following way: for any $x=\sum_{i=0}^\infty x_i p^i\in \zp$,
	$$
	Y_i(x)=\begin{cases}
		1 & \text{if } (x_{2i},x_{2i+1})=(0,0),\\
		0 & \text{if } (x_{2i},x_{2i+1})\neq (0,0).
	\end{cases}
	$$
	Notice that the random variables $(Y_i)_{i\in \mathbb N_0}$ are independent and identically distributed with $E[Y_i]=\frac{1}{p^2}$ for every $i\in \mathbb N_0$.
	Thus, by the strong law of large numbers, the set
	$$
	D=\left\{x=\sum_{i=0}^\infty x_i p^i\in \zp \colon \exists\lim_{n\rightarrow\infty} \frac{\sum_{i=0}^{n-1} Y_i(x)}{n}=\frac{1}{p^2} \right\}
	$$
	has measure $1$.
	Now, since for every $i\in \mathbb N_0$, $Y_i(x)=0$ for all $x=\sum_{j=0}^\infty x_j p^j$ that satisfy $(x_{2j},x_{2j+1})\neq (0,0)$ for each $j\in \mathbb N_0$, we have that $\lim_{n\to \infty} \frac{\sum_{i=0}^{n-1} Y_i(x)}{n}=0$ for all such $x$.
	Hence, it is clear that $E:=\left\{\sum_{i=0}^\infty  x_i p^i\in \zp\colon (x_{2i},x_{2i+1})\neq (0,0)  \text{ for all } i\in \mathbb N_0 \right\}\subset \zp\setminus D$.
	Moreover, by construction $\text{card}(E)=\mathfrak{c}$.
	Notice that it is not obvious that $E$ is a Borel set since any Haar measure $\mu$ on $(\zp,\mathcal B)$ is not complete.
	However, as $\zp\setminus D$ is a null set, we have that if $E$ were a Borel set, then $E$ would be a null set.
	Let us prove that $E$ is a Borel set.
	Consider the finite field of $p$ elements $\mathbb F_p$ endowed with an absolute value $|\cdot|_T$--the trivial absolute value.
	Then $\mathbb F_p$ is a discrete topological space, which implies that the product space $\mathbb F_p^2:=\mathbb F_p\times \mathbb F_p$ has the discrete topology. (Recall that the finite product of discrete topological spaces has the discrete topology.)
	For every $n\in \mathbb N_0$, let $\pi_n\colon \zp\to \mathbb F_p^2$ be defined as $\pi_n(x)=(x_{2n},x_{2n+1})$ for every $x=\sum_{i=0}^\infty x_i p^i\in \zp$.
	Let $n\in \mathbb N_0$, $x\in \zp$ and fix $\varepsilon>0$.
	Take an integer $m>n$.
	Then for every $y\in \zp$ such that $|x-y|_p<p^{-(2m+1)}$ we have that $\pi_n(x)=\pi_n(y)$, i.e., $|\pi_n(x)-\pi_n(y)|_T=0<\varepsilon$.
	Hence $\pi_n$ is continuous.
	Note that $E=\bigcap_{n=0}^\infty \pi_n^{-1}(\{(x,y)\in \mathbb F_p^2 \colon (x,y)\neq (0,0) \})$, where $\pi_n^{-1}(\{(x,y)\in \mathbb F_p^2 \colon (x,y)\neq (0,0) \})$ is closed since $\pi_n$ is continuous and $\{(x,y)\in \mathbb F_p^2 \colon (x,y)\neq (0,0) \}$ is closed in $\mathbb F_p^2$.
	Hence, $E$ is closed since it is the countable intersection of closed sets and, therefore, a Borel set.
	
	Let us analyze now the differentiability of $f$.
	On the one hand, if for $x=\sum_{i=0}^\infty x_i p^i\in \zp$ there exists $i_0\in\mathbb N_0$ such that $(x_{2i_0},x_{2i_0+1})=(0,0)$, then it is clear that $f$ is constant on some neighborhood of $x$, and hence differentiable at $x$.
	On the other hand, if $f$ were differentiable at $x=\sum_{i=0}^\infty x_i p^i\in \zp$ satisfying $(x_{2i},x_{2i+1})\neq (0,0)$ for all $i\in \mathbb N_0$, then we would have $f^\prime (x)=1$.
	Assume, by means of contradiction, that $f$ is differentiable at $x$.
	For every $n\in \mathbb N_0$, take $\overline{x}_n:=\sum_{i=0}^{2n+1} x_i p^i + \sum_{i=2n+4}^\infty y_i p^i$ with $y_{2n+4}\neq 0$, then
	\begin{align*}
		\frac{f(x)-f(\overline{x}_n)}{x-\overline{x}_n} & = \frac{\sum_{i=0}^\infty x_i p^i-\sum_{i=0}^{2n+1} x_i p^i}{\sum_{i=2n+2}^\infty x_{i} p^{i}-\sum_{i=2n+4}^\infty y_i p^{i} }\\
		& = \frac{\sum_{i=2n+2}^\infty x_{i} p^{i}}{\sum_{i=2n+2}^\infty x_{i} p^{i}-\sum_{i=2n+4}^\infty y_i p^{i}}\\
		& = \frac{\sum_{i=2n+2}^\infty x_{i} p^{i}-\sum_{i=2n+4}^\infty y_i p^{i}+\sum_{i=2n+4}^\infty y_i p^{i} }{\sum_{i=2n+2}^\infty x_{i} p^{i}-\sum_{i=2n+4}^\infty y_i p^{i}}\\
		& = 1+\frac{\sum_{i=2n+4}^\infty y_i p^{i}}{\sum_{i=2n+2}^\infty x_{i} p^{i}-\sum_{i=2n+4}^\infty y_i p^{i}}.
	\end{align*}
	Now, as $y_{2n+4} \neq 0$, we have
		$$
		\left|\frac{\sum_{i=2n+4}^\infty y_i p^{i}}{\sum_{i=2n+2}^\infty x_{i} p^{i}-\sum_{i=2n+4}^\infty y_i p^{i}} \right|_p=\begin{cases}
		p^{-2} & \text{if } x_{2n+2}\neq 0,\\
		p^{-1} & \text{if } x_{2n+3}\neq 0.
		\end{cases}
		$$
	Thus we have $\lim_{n\to \infty} |x-\overline{x}_n|_p=0$ and $\lim_{n\to \infty} \left|\frac{f(x)-f(\overline{x}_n)}{x-\overline{x}_n}-1\right|_p\geq p^{-2}\neq 0$, a contradiction.
	
	Let us define the function $g\colon \zp\rightarrow \qp$ by:
	$$
	g(x)=\begin{cases}
		p^n f(x^\prime) & \text{if } x=p^n+p^{n+1} x^\prime \text{ with } n\in\mathbb N \text{ and } x^\prime \in \zp,\\
		0 & \text{otherwise}.
	\end{cases}
	$$
	Notice that $g$ is well defined since the sets $B_n:=p^n+p^{n+1}\zp$ are pairwise disjoint.
	(The sets $B_n$ are the closed balls $\overline{B}_{p^{-(n+1)}}(p^n)$.)
	If $x\in \zp\setminus \left(\{0\} \cup \bigcup_{n=1}^\infty B_n\right)$, then there exists an open neighborhood $U^x$ of $x$ such that $g$ is identically zero on $U^x$, i.e., $g$ is differentiable at $x$.
	Now, as $g(p^n+p^{n+1}x)=p^n f(x)$ for every $n\in\mathbb N$ and $x\in \zp$, and since $f$ is continuous, it is obvious that $g$ is continuous on $\bigcup_{n=1}^\infty B_n$.
	Moreover, $g$ is also continuous at $0$.
	To prove it fix $\varepsilon>0$ and take $n\in \mathbb N$ such that $p^{-n}<\varepsilon$.
	If $x\in \zp$ is such that $|x|_p=p^{-n}$, then $x=x_n p^n+p^{n+1}x^\prime$ with $x_n\neq 0$.
	Furthermore,
		\begin{align*}
			|g(0)-g(x)|_p & =\begin{cases}
			|p^n f(x^\prime)|_p & \text{if } x_n=1,\\
			0 & \text{otherwise},
			\end{cases}
			=\begin{cases}
			p^{-n}|f(x^\prime)|_p & \text{if } x_n=1,\\
			0 & \text{otherwise}.
			\end{cases}
		\end{align*}
	Hence, $|g(0)-g(x)|_p\leq p^{-n}<\varepsilon$.
	Therefore we have proven that $g$ is continuous on $\zp$.
	Moreover, $g$ is differentiable also on $\bigcup_{n=1}^\infty (B_n\setminus E_n)$ (with bounded derivative) as $f$ is differentiable on $\bigcup_{n=1}^\infty (B_n\setminus E_n)$, where $E_n:=p^n+p^{n+1}E$; and $g$ is not differentiable on $\bigcup_{n=1}^\infty E_n$ since $f$ is not differentiable on $\bigcup_{n=1}^\infty E_n$.
	Notice that once again $\text{card}(E_n)=\mathfrak{c}$ for every $n\in \mathbb N_0$.
	
	Let us prove that $E_n$ is a Borel set with $\mu(E_n)=0$ for every $n\in \mathbb N$.
	To do so, let us consider the restricted measure $\mu_n=p^{n+1}\mu$ on the measurable space $(B_n,\mathcal B_n)$, where $\mathcal B_n$ is the $\sigma$-algebra of all Borel subsets of $B_n$.
	Notice that $\mathcal B_n=\{B\cap B_n\colon B\in \mathcal B \}$ and $(B_n,\mathcal B_n,\mu_n)$ is a probability space.
	Define now for every $i\in \mathbb N_0$ the random variables $Y_{n,i}\colon B_n\rightarrow\{0,1\}$ as follows: for $x=p^n+p^{n+1}\sum_{i=0}^\infty x_i p^i\in B_n$, we have
	$$
	Y_{n,i}(x)=\begin{cases}
		1 & \text{if } (x_{2i},x_{2i+1})=(0,0),\\
		0 & \text{if } (x_{2i},x_{2i+1}) \neq (0,0).
	\end{cases}
	$$
	Once again the random variables $(Y_{n,i})_{i\in \mathbb N_0}$ are independent and identically distributed with $E[Y_{n,i}]=\frac{1}{p^2}$ for every $i\in\mathbb N_0$.
	Thus, the set 
	$$
	\left\{x=p^n+p^{n+1}\sum_{i=0}^\infty x_i p^i\in B_n \colon \exists\lim_{m\rightarrow\infty} \frac{\sum_{i=0}^{m-1} Y_{n,i}(x)}{m}=\frac{1}{p} \right\}=p^n+p^{n+1}D=:D_n
	$$
	is a full set for $\mu_n$.
	By considering for each $k\in \mathbb N_0$ the function $\pi_{n,k}\colon B_n\to \mathbb F_p^2$ given by $\pi_{n,k}(x)=(x_{2k},x_{2k+1})$ for every $x=p^n+p^{n+1}\sum_{i=0}^\infty x_i p^i\in B_n$ and applying similar arguments as above, we have that $\pi_{n,k}$ is continuous.
	Hence $E_n=\bigcap_{k=0}^\infty \pi_{n,k}(\{(x,y)\in \mathbb F_p^2\colon (x,y) \neq (0,0) \})$ is once again a Borel set.
	Furthermore, since $E_n\subset B_n\setminus D_n$ we have that $E_n$ is a null set for $\mu_n$.
	Thus $\mu(E_n)=p^{-(n+1)}\mu_n(E_n)=0$ for every $n\in \mathbb N$.
	
	Finally let us prove that $g$ is not differentiable at $0$.
	Since every neighborhood containing $0$ on $\zp$ contains points $x$ such that $g(x)=0$, if $g$ were differentiable at $0$ then $g^\prime (0)=0$.
	Assume that $g$ is differentiable at $0$.
	As $p^n+p^{n+1}\sum_{i=0}^\infty p^i=p^n \sum_{i=0}^\infty p^i\in B_n$ for every $n\in \mathbb N$, we have that
		$$
		\left|\frac{g\left(p^n+p^{n+1}\sum_{i=0}^\infty p^i \right)-g(0)}{p^n+p^{n+1}\sum_{i=0}^\infty p^i} \right|_p=\left|\frac{p^n f\left(\sum_{i=0}^\infty p^i \right)}{p^n \sum_{i=0}^\infty p^i} \right|_p=\left|\frac{p^n \sum_{i=0}^\infty p^i}{p^n \sum_{i=0}^\infty p^i}  \right|_p=1,
		$$
	where $\lim_{n\to \infty } \left|p^n+p^{n+1}\sum_{i=0}^\infty p^i \right|_p=\lim_{n\to \infty} p^{-n}=0$, a contradiction.

	For every $N\in \mathcal N$, let us define $f_N\colon \zp\rightarrow \qp$ by:
	$$
	f_N(x)=g(x)\sum_{n\in N} 1_{B_n}(x).
	$$
	Since $B_n\cap B_m=\emptyset$ for every distinct $n,m\in \mathbb N$, the function $f_N$ is well defined.
	Furthermore, since each $N\in \mathcal N$ is infinite, we can apply the above arguments to prove that $f_N$ is continuous and differentiable on a full set for $\mu$ with bounded derivative but not differentiable on the complement having cardinality $\mathfrak c$.
	
	It remains to prove that the functions in $V=\{f_N\colon N\in \mathcal N \}$ are linearly independent over $\qp$ and any nonzero linear combination over $\qp$ of the functions in $V$ satisfies the necessary properties.
	Let $F:=\sum_{i=1}^k a_i f_{N_i}$, where $k\in \mathbb N$, $a_1,\ldots,a_k\in \qp$, and $N_1,\ldots,N_k\in \mathcal N$ are distinct.
	We begin by showing the linear independence.
	Assume that $F\equiv 0$.
	Fix $n\in N_1^1\cap N_2^0\cap \cdots \cap N_k^0$ and take $x=p^n+p^{n+1}\sum_{i=0}^\infty p^i\in B_n$.
	Then, $0=F(x)=a_1 f_{N_1}(x)=a_1 \sum_{i=0}^\infty p^i$ if and only if $a_1=0$.
	By repeating the same argument, it is easy to see that $a_i=0$ for every $i\in \{1,\ldots,k \}$.
	Assume now that $a_i\neq 0$ for every $i\in \{1,\ldots,k \}$.
	Then $F$ is continuous but also differentiable on
		$$
		\Delta:=\zp\setminus \left(\{0\}\cup \left(\bigcup_{n\in \bigcup_{i=1}^k N_i} E_n \right) \right)
		$$
	(with bounded derivative).
		Applying similar arguments as above, we have that $F$ is not differentiable at $0$.
	Let $x\in E_n$ with $n\in \bigcup_{i=1}^k N_i$.
	We will analyze the differentiability of $F$ at $x$ depending on the values that $F$ takes on $B_n$.
	We have two possible cases.
	
	\textit{Case 1:} If $F$ is identically $0$ on $B_n$, then $F$ is differentiable at $x$.
	
	\textit{Case 2:} If $F$ is not identically $0$ on $B_n$, then there exists $a\in \qp\setminus\{0\}$ such that $F=ag$.
	Hence $F$ is not differentiable at $x$ since $g$ is not differentiable at $x$.
	Notice that Case 2 is always satisfied.
	
	To finish the proof, it is enough to show that $\Delta$ is a full set for $\mu$, but this is an immediate consequence of the fact that $\{0\}\cup \left(\bigcup_{n\in \bigcup_{i=1}^k N_i} E_n \right)$ is the countable union of null sets for $\mu$ since it implies that
		$$
		\mu\left(\{0\}\cup \left(\bigcup_{n\in \bigcup_{i=1}^k N_i} E_n \right) \right)=0.
		$$
\end{proof}


\begin{bibdiv}
\begin{biblist}
	
\bib{abms}{article}{
	AUTHOR = {Ara\'{u}jo, G.},
	AUTHOR = {Bernal-Gonz\'{a}lez, L.},
	AUTHOR = {Mu\~{n}oz-Fern\'{a}ndez, G. A.},
	AUTHOR = {Seoane-Sep\'{u}lveda, J. B.},
	TITLE = {Lineability in sequence and function spaces},
	JOURNAL = {Studia Math.},
	VOLUME = {237},
	YEAR = {2017},
	NUMBER = {2},
	PAGES = {119--136},
}
	
	

\bib{ar}{book}{
	author={Aron, R.M.},
	author={Bernal Gonz\'{a}lez, L.},
	author={Pellegrino, D.M.},
	author={Seoane Sep\'{u}lveda, J.B.},
	title={Lineability: the search for linearity in mathematics},
	series={Monographs and Research Notes in Mathematics},
	publisher={CRC Press, Boca Raton, FL},
	date={2016},
	pages={xix+308},
	isbn={978-1-4822-9909-0},
}



\bib{AGS2005}{article}{
	author={Aron, R. M.},
	author={Gurariy, V. I.},
	author={Seoane-Sep\'{u}lveda, J. B.},
	title={Lineability and spaceability of sets of functions on $\mathbb R$},
	journal={Proc. Amer. Math. Soc.},
	volume={133},
	date={2005},
	number={3},
	pages={795--803},
}

\bib{ar2}{article}{
	author={Aron, R.M.},
	author={P\'{e}rez-Garc\'{\i }a, D.},
	author={Seoane-Sep\'{u}lveda, J.B.},
	title={Algebrability of the set of non-convergent Fourier series},
	journal={Studia Math.},
	volume={175},
	date={2006},
	number={1},
	pages={83--90},
}

\bib{bbf}{article}{
	AUTHOR = {Balcerzak, Marek},
		AUTHOR = {Bartoszewicz, Artur},
		AUTHOR = {Filipczak, Ma\l gorzata},
	TITLE = {Nonseparable spaceability and strong algebrability of sets of
		continuous singular functions},
	JOURNAL = {J. Math. Anal. Appl.},
	VOLUME = {407},
	YEAR = {2013},
	NUMBER = {2},
	PAGES = {263--269},
	}
\bib{bbfg}{article}{
	AUTHOR = {Bartoszewicz, Artur},
		AUTHOR = {Bienias, Marek},
		AUTHOR = {Filipczak, Ma\l gorzata},
		AUTHOR = {G\l \c{a}b, Szymon},
	TITLE = {Strong {$\germ{c}$}-algebrability of strong
		{S}ierpi\'{n}ski-{Z}ygmund, smooth nowhere analytic and other sets
		of functions},
	JOURNAL = {J. Math. Anal. Appl.},
	VOLUME = {412},
	YEAR = {2014},
	NUMBER = {2},
	PAGES = {620--630},
}


\bib{barglab}{article}{
	author={Bartoszewicz, A.},
	AUTHOR = {Bienias, M.},
	author={G\l \c {a}b, S.},
	TITLE = {Independent {B}ernstein sets and algebraic constructions},
	JOURNAL = {J. Math. Anal. Appl.},
	VOLUME = {393},
	YEAR = {2012},
	NUMBER = {1},
	PAGES = {138--143},
}



\bib{bar20}{article}{
	AUTHOR = {Bartoszewicz, Artur},
	AUTHOR ={ Filipczak, Ma\l gorzata},
	AUTHOR = {Terepeta, Ma\l gorzata},
	TITLE = {Lineability of {L}inearly {S}ensitive {F}unctions},
	JOURNAL = {Results Math.},
	VOLUME = {75},
	YEAR = {2020},
	NUMBER = {2},
	PAGES = {Paper No. 64},
}


\bib{ba3}{article}{
	author={Bartoszewicz, A.},
	author={G\l \c {a}b, S.},
	title={Strong algebrability of sets of sequences and functions},
	journal={Proc. Amer. Math. Soc.},
	volume={141},
	date={2013},
	number={3},
	pages={827--835},
}





\bib{bay}{article}{
	AUTHOR = {Bayart, Fr\'{e}d\'{e}ric},
	TITLE = {Linearity of sets of strange functions},
	JOURNAL = {Michigan Math. J.},
	VOLUME = {53},
	YEAR = {2005},
	NUMBER = {2},
	PAGES = {291--303},
	}


\bib{bq}{article}{
	AUTHOR = {Bayart, Fr\'{e}d\'{e}ric},
	AUTHOR = {Quarta, Lucas},
	TITLE = {Algebras in sets of queer functions},
	JOURNAL = {Israel J. Math.},
	VOLUME = {158},
	YEAR = {2007},
	PAGES = {285--296},
}

\bib{bbls}{article}{
	AUTHOR = {Bernal-Gonz\'{a}lez, L.},
	AUTHOR = {Bonilla, A.},
	AUTHOR = {L\'{o}pez-Salazar, J.},
	AUTHOR = {Seoane-Sep\'{u}lveda, J. B.},
	TITLE = {Nowhere h\"{o}lderian functions and {P}ringsheim singular
		functions in the disc algebra},
	JOURNAL = {Monatsh. Math.},
	VOLUME = {188},
	YEAR = {2019},
	NUMBER = {4},
	PAGES = {591--609},
}


\bib{TAMS2020}{article}{
	author={Bernal-Gonz\'{a}lez, L.},
	author={Cabana-M\'{e}ndez, H. J.},
	author={Mu\~{n}oz-Fern\'{a}ndez, G. A.},
	author={Seoane-Sep\'{u}lveda, J. B.},
	title={On the dimension of subspaces of continuous functions attaining
		their maximum finitely many times},
	journal={Trans. Amer. Math. Soc.},
	volume={373},
	date={2020},
	number={5},
	pages={3063--3083},
}

\bib{b20}{article}{
	author={Bernal-Gonz\'{a}lez, L.},
	author={Mu\~{n}oz-Fern\'{a}ndez, G. A.},
	author={Rodr\'{\i }guez-Vidanes, D. L.},
	author={Seoane-Sep\'{u}lveda, J. B.},
	title={Algebraic genericity within the class of sup-measurable functions},
	journal={J. Math. Anal. Appl.},
	volume={483},
	date={2020},
	number={1},
	pages={123--576},
}

\bib{bo}{article}{
	AUTHOR = {Bernal-Gonz\'{a}lez, Luis},
	AUTHOR = {Ord\'{o}\~{n}ez Cabrera, Manuel},
	TITLE = {Lineability criteria, with applications},
	JOURNAL = {J. Funct. Anal.},
	VOLUME = {266},
	YEAR = {2014},
	NUMBER = {6},
	PAGES = {3997--4025},
}

\bib{survey}{article}{
	author={Bernal-Gonz\'{a}lez, L.},
	author={Pellegrino, D.},
	author={Seoane-Sep\'{u}lveda, J.B.},
	title={Linear subsets of nonlinear sets in topological vector spaces},
	journal={Bull. Amer. Math. Soc. (N.S.)},
	volume={51},
	date={2014},
	number={1},
	pages={71--130},
}





\bib{bns}{article}{
	AUTHOR = {Biehler, N.},
	AUTHOR = {Nestoridis, V.},
	AUTHOR = {Stavrianidi, A.},
	TITLE = {Algebraic genericity of frequently universal harmonic
		functions on trees},
	JOURNAL = {J. Math. Anal. Appl.},
	VOLUME = {489},
	YEAR = {2020},
	NUMBER = {1},
	PAGES = {124132, 11},
}

\bib{Bi}{book}{
	author={Billingsley, Patrick},
	title={Probability and measure},
	series={Wiley Series in Probability and Mathematical Statistics},
	edition={3},
	note={A Wiley-Interscience Publication},
	publisher={John Wiley \& Sons, Inc., New York},
	date={1995},
	pages={xiv+593},
}

\bib{cgp}{article}{
	AUTHOR = {Calder\'{o}n-Moreno, M. C.},
	AUTHOR = {Gerlach-Mena, P. J.},
	AUTHOR = {Prado-Bassas, J. A.},
	TITLE = {Lineability and modes of convergence},
	JOURNAL = {Rev. R. Acad. Cienc. Exactas F\'{\i}s. Nat. Ser. A Mat. RACSAM},
	VOLUME = {114},
	YEAR = {2020},
	NUMBER = {1},
	PAGES = {Paper No. 18, 12},
}



\bib{CS2019}{article}{
	author={Ciesielski, Krzysztof C.},
	author={Seoane-Sep\'{u}lveda, Juan B.},
	title={Differentiability versus continuity: restriction and extension theorems and monstrous examples},
	journal={Bull. Amer. Math. Soc. (N.S.)},
	volume={56},
	date={2019},
	number={2},
	pages={211--260},
}

\bib{egs}{article}{
	AUTHOR = {Enflo, Per H.},
		AUTHOR = {Gurariy, Vladimir I.},
		AUTHOR = {Seoane-Sep\'{u}lveda,	Juan B.},
	TITLE = {Some results and open questions on spaceability in function
		spaces},
	JOURNAL = {Trans. Amer. Math. Soc.},
	VOLUME = {366},
	YEAR = {2014},
	NUMBER = {2},
	PAGES = {611--625},
}




\bib{fmrs}{article}{
	AUTHOR = {Fern\'{a}ndez-S\'{a}nchez, J.},  
	author={Maghsoudi, S.},
	author={Rodr\'{\i}guez-Vidanes, D. L.},
	author={Seoane-Sep{\'u}lveda, J. B.},
	TITLE = {Classical vs. non-Archimedean analysis. An approach via algebraic genericity},
   status={Preprint (2021)},
}

\bib{preprint-June-2019}{article}{
	author={Fern\'andez-S\'anchez, J.},
	author={Mart\'inez-G\'omez, M.E.},
	author={Seoane-Sep\'ulveda, J.B.},
	title={Algebraic genericity and special properties within sequence spaces and series},
	journal={Preprint (2019)},
}

\bib{fich}{article}{
	author={Fichtenholz, G.},
	author={Kantorovich, L.},
	title={Sur les op\'{e}rations dans l'espace des functions born\'{e}es},
	journal={Studia Math.},
	volume={5},
	date={1934},
	number={},
	pages={69-98},
	}

\bib{Fo}{book}{
	author={Folland, Gerald B.},
	title={A course in abstract harmonic analysis},
	series={Textbooks in Mathematics},
	edition={2},
	publisher={CRC Press, Boca Raton, FL},
	date={2016},
	pages={xiii+305 pp.+loose errata},
}

\bib{g7}{article}{
	AUTHOR = {Garc\'{\i}a-Pacheco, F. J.},
	AUTHOR = {Palmberg, N.},
	AUTHOR = {Seoane-Sep\'{u}lveda, J. B.},
	TITLE = {Lineability and algebrability of pathological phenomena in analysis},
	JOURNAL = {J. Math. Anal. Appl.},
	VOLUME = {326},
	YEAR = {2007},
	NUMBER = {2},
	PAGES = {929--939},
}



\bib{go}{book}{
	author={Gouv\^{e}a, F.Q.},
	title={$p$-adic numbers, An introduction},
	series={Universitext},
	edition={2},
	publisher={Springer-Verlag, Berlin},
	date={1997},
	pages={vi+298},
}



\bib{gur1}{article}{
	author={Gurari\u {\i }, V.I.},
	title={Subspaces and bases in spaces of continuous functions},
	language={Russian},
	journal={Dokl. Akad. Nauk SSSR},
	volume={167},
	date={1966},
	pages={971--973},
}

\bib{h}{article}{
	author={Hausdorff, F.},
	title={Uber zwei Satze von G. Fichtenholz und L. Kantorovich},
	language={German},
	journal={Studia Math.},
	volume={6},
	date={1936},
	pages={18--19},
}

\bib{jms}{article}{
	AUTHOR = {Jim\'{e}nez-Rodr\'{\i}guez, P.},
		AUTHOR = {Mu\~{n}oz-Fern\'{a}ndez, G. A.},
		AUTHOR = {Seoane-Sep\'{u}lveda, J. B.},
			TITLE = {Non-{L}ipschitz functions with bounded gradient and related
		problems},
	JOURNAL = {Linear Algebra Appl.},
	VOLUME = {437},
	YEAR = {2012},
	NUMBER = {4},
	PAGES = {1174--1181},
	}



\bib{kato}{book}{
	AUTHOR = {Katok, S.},
	TITLE = {{$p$}-adic analysis compared with real},
	SERIES = {Student Mathematical Library},
	VOLUME = {37},
	PUBLISHER = {American Mathematical Society, Providence, RI; Mathematics
		Advanced Study Semesters, University Park, PA},
	YEAR = {2007},
	PAGES = {xiv+152},
}

\bib{jms1}{article}{
	author={Khodabendehlou, J.},
	author={Maghsoudi, S.},
	author={Seoane-Sep{\'u}lveda, J.B.},
	title={Algebraic genericity and summability within the non-Archimedean setting},
	journal={Rev. R. Acad. Cienc. Exactas F\'{\i}s. Nat. Ser. A Mat. RACSAM},
	volume={115},
	number={21},
	date={2021},
}



\bib{jms2}{article}{
	author={Khodabendehlou, J.},
	author={Maghsoudi, S.},
	author={Seoane-Sep{\'u}lveda, J.B.},
	title={Lineability and algebrability within $p$-adic function spaces},
	journal={Bull. Belg. Math. Soc. Simon Stevin},
	volume={27},
	pages={711–729},
	date={2020},
}

\bib{jms3}{article}{
	author={Khodabendehlou, J.},
	author={Maghsoudi, S.},
	author={Seoane-Sep{\'u}lveda, J.B.},
	title={Lineability, continuity, and antiderivatives in the
		non-Archimedean setting},
	journal={Canad. Math. Bull.},
	volume={64},
	date={2021},
	number={3},
	pages={638--650},
}


\bib{levinemilman1940}{article}{
	author={Levine, B.},
	author={Milman, D.},
	title={On linear sets in space $C$ consisting of functions of bounded variation},
	language={Russian, with English summary},
	journal={Comm. Inst. Sci. Math. M\'{e}c. Univ. Kharkoff [Zapiski Inst. Mat. Mech.] (4)},
	volume={16},
	date={1940},
	pages={102--105},
}




\bib{mahler}{book}{
	AUTHOR = {Mahler, K.},
	TITLE = {$p$-adic numbers and their functions},
	SERIES = {Cambridge Tracts in Mathematics},
	VOLUME = {76},
	EDITION = {Second},
	PUBLISHER = {Cambridge University Press, Cambridge-New York},
	YEAR = {1981},
	PAGES = {xi+320},
	}

\bib{Ma2}{article}{
	author={Mahler, K.},
	title={An interpolation series for continuous functions of a $p$-adic
		variable},
	journal={J. Reine Angew. Math.},
	volume={199},
	date={1958},
	pages={23--34},
}

\bib{ms}{article}{
	AUTHOR = {Moothathu, T. K. Subrahmonian},
	TITLE = {Lineability in the sets of {B}aire and continuous real
		functions},
	JOURNAL = {Topology Appl.},
	VOLUME = {235},
	YEAR = {2018},
	PAGES = {83--91},
}



\bib{nak}{article}{
	AUTHOR = {Natkaniec, Tomasz},
	TITLE = {On lineability of families of non-measurable functions of two
		variable},
	JOURNAL = {Rev. R. Acad. Cienc. Exactas F\'{\i}s. Nat. Ser. A Mat. RACSAM},
	VOLUME = {115},
	YEAR = {2021},
	NUMBER = {1},
	PAGES = {Paper No. 33, 10},
}



\bib{ro}{book}{
	AUTHOR = {Robert, A. M.},
	TITLE = {A course in $p$-adic analysis},
	SERIES = {Graduate Texts in Mathematics},
	VOLUME = {198},
	PUBLISHER = {Springer-Verlag, New York},
	YEAR = {2000},
	PAGES = {xvi+437},
}

\bib{sc}{book}{
	author={Schikhof, W.H.},
	title={Ultrametric calculus, An introduction to $p$-adic analysis},
	series={Cambridge Studies in Advanced Mathematics},
	volume={4},
	publisher={Cambridge University Press, Cambridge},
	date={1984},
	pages={viii+306},
}



\bib{juanksu}{book}{
	author={Seoane-Sep\'{u}lveda, J.B.},
	title={Chaos and lineability of pathological phenomena in analysis},
	note={Thesis (Ph.D.)--Kent State University},
	publisher={ProQuest LLC, Ann Arbor, MI},
	date={2006},
	pages={139},
}




\bib{va}{book}{
	author={van Rooij, A.C.M.},
	title={Non-Archimedean functional analysis},
	series={Monographs and Textbooks in Pure and Applied Math.},
	volume={51},
	publisher={Marcel Dekker, Inc., New York},
	date={1978},
	pages={x+404},
}
\end{biblist}
\end{bibdiv}
\end{document}